\documentclass{article}
\usepackage{amssymb}
\usepackage{amsthm}
\usepackage{eufrak}
\usepackage[all]{xy}
\usepackage{makeidx}

\makeindex
\title{LATTICES AND COHOMOLOGY.}
\author{Luis Arenas-Carmona. 
\footnote{Supported by Fondecyt, proyecto No. 3010018,
and the Chilean Catedra Presidencial in Number Theory.}\\
\\ 
\small\textit{Universidad de Chile,}\\
\small\textit{Facultad de Ciencias.}\\
\small{Casilla 653, Santiago,Chile.}\\
\small{learenas@uchile.cl}}
\date{}

\theoremstyle{plain}
\newtheorem{prop}{Proposition}[section]
\newtheorem{propi}{Proposition}[section]
\newtheorem{lem}[prop]{Lemma}
\theoremstyle{definition}
\newtheorem{ex}[prop]{Example}
\newtheorem{dfn}[prop]{Definition}
\newtheorem{rmk}[prop]{Remark}

\newcommand\End{\textnormal{End}}

\newcommand\Hom{\textnormal{Hom}}
\newcommand\smallgen{\textnormal{\tiny{gen}}}

\newcommand\im{\textnormal{im}}

\newcommand\ckgen[3]{C_{\smallgen}(#1,#2,#3)}

\newcommand\ad{\mathbb A}

\newcommand\la{\Lambda}
\newcommand\bark{\bar{k}}

\newcommand\uno{\{1\}}

\newcommand\oink{\mathcal O}

\newcommand\Q{\mathbb Q}

\newcommand\enteri{\mathbb Z}

\newcommand\gal{\mathcal G}

\newcommand\galbark{{\mathcal G_{\bark/k}}}

\newcommand\galbarkk{\galbark}
\newcommand\galkk{{\mathcal G_{K/k}}}
\newcommand\galkvkv{{\mathcal G_{K_w/k_v}}}
\newcommand\galv{{\mathcal G_w}}
\newcommand\tn{\mathfrak{T}}
\newcommand\dete{\textnormal{det}}

\newcommand\vaa{\longrightarrow}

\newcommand\grupi{\mathbb G}

\newcommand\ideala{\mathcal A}
\newcommand\idealp{\wp}
\newcommand\idealP{\mathcal P}

\newcommand\orth[3]{\oink_{#1,#2}(#3)}

\newcommand\orthi[4]{\oink_{#1,#2}^{#4}(#3)}

\newcommand\orthomega{\oink_n(Q)}

\newcommand\orthkbar{\orth n{\bark}Q}

\newcommand\linei[4]{GL_{#2}^{#4}(#3)}
\newcommand\linear[3]{\linei {#1}{#2}{#3}{}}
\newcommand\linek{\linear nkV}
\newcommand\lineE{\linear nEV}
\newcommand\lineK{\linear nKV}

\newcommand\lineklam{\linei nkV\la}
\newcommand\lineKlam{\linei nKV\la}

\newcommand\lineKvlam{\linei n{K_w}V\la}
\newcommand\lineomega{GL(V)}
\newcommand\linekbarlam{\linei n{\bark}V\la}

\newcommand\speciali[3]{SL_{#1}^{#3}(#2)}

\newcommand\specialKlam{\speciali KV\la}

\newcommand\specialomega{SL(V)}

\newcommand\tensor{\bigotimes_{\oink_k}}

\newcommand\lk{\mathcal{L}_{\textnormal{\tiny def}}}
\newcommand\lfree{\mathcal{L}_{\textnormal{\tiny fr}}}
\newcommand\lvk{\mathcal{L}^V}
\newcommand\lkvk{\mathcal{L}_{\textnormal{\tiny def}}^V}
\newcommand\lfreevk{\mathcal{L}_{\textnormal{\tiny fr}}^V}

\newcommand\Lk[3]{\mathcal{L}_{\textnormal{\tiny def}}(#1,#2,#3)}
\newcommand\Lfree[3]{\mathcal{L}_{\textnormal{\tiny fr}}(#1,#2,#3)}
\newcommand\Lvk[3]{\mathcal{L}^V(#1,#2,#3)}
\newcommand\Lkvk[3]{\mathcal{L}_{\textnormal{\tiny def}}^V(#1,#2,#3)}
\newcommand\Lfreevk[3]{\mathcal{L}_{\textnormal{\tiny fr}}^V
(#1,#2,#3)}

\newcommand\ksobrekn{k^*/(k^*)^n}

\newcommand\gtilde{\widetilde{G}}

\newcommand\gggr[2]{G_{#1}^{#2}}

\newcommand\gkla{\gggr k\la}
\newcommand\gakla{\gggr {\ad_k}\la}
\newcommand\gKla{\gggr K\la}

\newcommand\gKlav{\gggr {K_w}\la}

\newcommand\gK{G_K}
\newcommand\gk{G_k}

\newcommand\gkv{G_{k_v}}
\newcommand\gak{G_{\ad_k}}

\begin{document}
\maketitle

\begin{abstract}
We give an interpretation of the cohomology of an arithmetically defined
group  as a set of equivalence classes of lattices. We  use this 
interpretation to give a simpler proof of the connection established 
in \cite{Rohlfs} between genus and cohomology.
\end{abstract}


\section{Introduction.}

Galois cohomology is a fundamental tool for the 
classification of certain algebraic structures.

To be precise, let $k$ be a field, $G$ a linear algebraic group
acting on a space $V$, both defined over $k$.
It is known (see \cite{kneser2}), that if $G$ is
 defined as
 the set of automorphisms of a tensor $\tau$ on $V$,
 e.g., a quadratic form or an algebra structure, 
 the cohomology set $H^1(K/k,G_K)$ classifies the
 $K/k$-forms of $\tau$,
 i.e., those tensors of the same
 type also defined over $k$ that 
 become isomorphic to $\tau$
over the larger field $K$ (see section \ref{kkforms}).

It would be reasonable to expect, therefore, that a similar theory
were available for structures for which the corresponding automorphism
group is not an algebraic group but an 
arithmetically defined subgroup of an algebraic group.

Such a theory is already hinted at in \cite{Rohlfs}.  In this
reference, two finiteness results are proven. The first one deals with
the finiteness of the local cohomology set $H^1(\gal_w,\Gamma_w)$,
for an arithmetically defined group $\Gamma$.
Notations are as in \cite{Rohlfs}.
The second one deals with the finiteness of the kernel of the map
\[
H^1(\gal,\Gamma)\vaa \prod_{v\textnormal{ place of
}k}H(\gal_{w(v)},\Gamma_{w(v)}),
\]
where we have fixed a place $w(v)$ 
of $K$ dividing each place $v$ of $k$.  
It is the proof of the second result which requires
expressing the given kernel in terms of the set of double cosets
\[
G_k\backslash G_{\ad_k}/\prod_w\Gamma_w
\]
 (see corollary 3.3 in \cite{Rohlfs}).
These double cosets  are the same ones that
classify the classes of lattices in a genus. This is the relation
we want to pursue.

In this paper, we show that the results in \cite{Rohlfs} are part of a 
much more general theory that relates cohomology sets of 
arithmetically defined groups  with certain equivalence 
classes of lattices in $V$.

Let $k$ be a local or number fields, $K/k$ a Galois extension 
with Galois group $\gal=\galkk$.
Let $V_K,G_K$ denote the sets of $K$-points of $V$ and $G$
(see section \ref{notations}).
We establish the following result:

\begin{propi}
There exists a correspondence between 
$G_k$-orbits of $\gal$-stable lattices in $V_K$,
that are in the same $G_K$-orbit,
and elements of 
\[
\ker\Big(H^1(\gal,G_K^{\la})\stackrel{i_*}{\longrightarrow}H^1(\gal,G_K)\Big),
\]
where $G_K^\la$ is the stabilizer of one particular lattice,
and $i_*$ the map induced by the inclusion
(see prop. \ref{result1}).
\end{propi}

The cocycles described above can be thought as equivalence classes of
\emph{lattices in the same
space}. We develop a concept of \emph{ lattices
in other spaces}, at least for the case that $G$ is the stabilizer 
of a family of tensors (see \ref{kkforms}).
In this context, the following result is obtained:

\begin{propi}
Assume $G$ is the 
stabilizer of a family of tensors $\tn$ on $V$.
Then, the set $H^1(\gal,\gKla)$ is 
in one-to-one correspondence 
with the set of $G_k$-orbits of $\gal$-invariant 
lattices in other spaces that are isomorphic over $K$. The set
\[
\ker\Big(
H^1(\gal,\gKla)\stackrel{i_*}{\vaa}H^1(\gal,\gK)\Big),
\]
where $i$ is the inclusion, corresponds to those lattices 
that are in the
same space as $\la_k$  (see prop. \ref{section3.4}).
\end{propi}

 The relation between the set of cohomology classes and the set of lattice
classes is established by taking
 advantage of the long exact sequence in cohomology that
 arises from a short exact sequence over $K$.

 We analyze the cohomology of the general linear group 
 and its relation to classification of lattices.
 In particular, one obtains that the set of $G_k$-orbits
 of free lattices, that are isomorphic over $K$ to a 
given lattice, corresponds to the set
\[
\ker\Bigg(H^1(\galkk,\gKla)\vaa
H^1\Big(\galkk,\lineKlam\Big)\Bigg)
\]
(see prop. \ref{proposition7}).

A lattice is defined over $k$, if it is generated by its $k$-points.
This is a stronger condition that $\galkk$-invariance.
 It is necessary to look at localizations 
to obtain a cohomological characterization
 of the set of lattice classes defined over $k$. 
We show that this set corresponds to the kernel
\[
\ker\Bigg(H^1(\galkk,\gKla)\vaa \prod_{v}
H^1\Big(\galv,\lineKvlam\Big)\Bigg),
\]
where $\galv$ is the local Galois group
(see prop. \ref{lkernel}).

 Section \ref{chapter6} is devoted to the study of the relation between
 integral cohomology and genera, ie., how
cohomology can be used  to study 
sets of lattices that become isometric over some extension.
This is expressed in terms of the notion of 
\emph{cohomological genus}.
In particular, we recover the results in \cite{Rohlfs}.

\section{Notations.}\label{notations}

In all of this article, $k,K,E$ denote number or local fields of
characteristic $0$, or algebraic extensions of them.
If $k$ is a number field, $\Pi(k)$ denotes the set of
places of $k$.
\begin{rmk}\label{omegagroups}
By an algebraic group, we mean a linear algebraic
group. All algebraic groups are assumed to be subgroups
of the general linear group of a vector space $V$, 
of finite dimension, over
a sufficiently large
algebraically closed
field $\Omega$ of characteristic $0$.
We assume that all localizations of number fields
inject into $\Omega$.
$G$ denotes an algebraic group over $\Omega$.
 $\lineomega,\specialomega$ denote the general and 
special linear groups over $\Omega$.
When we work over a
fixed local or number field $k$,
we say that $G$ is defined over $k$ if the equations defining $G$
have coefficients in $k$ (see section 2.1.1 in \cite{Pla}).
 \emph{This is the case for all groups considered
here}.
For any field $E$, $k\subseteq E\subseteq\Omega$, we write $G_E$ for the set
of $E$-points of $G$, e.g., if $G=\lineomega$, the set of $E$ points
is denoted $\lineE$. The same conventions apply to spaces and algebras.
All spaces and algebras are assumed to be finite dimensional.

Exceptions to this rule are the multiplicative and additive groups.
We denote the multiplicative group $\Omega^*$ by $\grupi_m$,
and the additive group $\Omega$ by $\grupi_a$, when considered as
algebraic groups. For the set of $k$-points
we write $k^*$,$k$. Instead of $(\grupi_m)_k,(\grupi_a)_k$.
\end{rmk}

The orthogonal group of a quadratic form $Q$ on $V$ is written
$\oink_{n}(Q)$ or $\oink_{n}(Q,V)$, where $n=\dim_\Omega(V)$.
The set of $E$-points is denoted $\oink_{n,E}(Q)$.

The field on which a particular lattice is defined is always written as
a subindex. If $K/k$ an extension of local or number fields
and $\la_k$ is a lattice in $V_k$, $\la_K$ denotes the $\oink_K$-lattice
in $V_K$ generated by $\la_k$.

If $G$ is an algebraic group  acting on a space $V$, both 
defined over $k$,
and $\la_k$ is a $\oink_k$-lattice on $V_k$, the stabilizer of $\la_k$
in $G_k$ is denoted $\gkla$. If $G=\lineomega$, this set is denoted
$\lineklam$. Similar conventions apply to 
special linear or orthogonal groups.

\begin{rmk}\label{fixw}
Whenever $K/k$ is a Galois extension of a number field $k$, and
 $v$ a place of $k$,
$w$ denotes a place of $K$ dividing $v$. We assume that one fixed such $w$
 has been chosen for every $v$. This convention is also applied for
infinite extension, e.g., $K=\bark$.
\end{rmk}

\begin{rmk}\label{galoisgroups}
$\galkk$ denotes the Galois group of the extension $K/k$.
If there is no risk of confusion, we write simply $\gal$.
If $K$ is not specified, we assume $K=\bar{k}$.
If $k$ is a number field and $v\in\Pi(k)$,
we also use the notation $\galv=\galkvkv$.
\end{rmk}

If $\Gamma$ is a group acting on a set $S$, $S/\Gamma$ denotes the set of orbits
and $S^\Gamma$ the set of invariant points. The action of $\gamma\in\Gamma$
is denoted $s\mapsto s^\gamma$, for $s\in S$.

\section{Review of Galois cohomology.}\label{ReviewDefinitions}

In this section we recall some of the cohomology results that we need.
The results in this section are found in chapter 1 in \cite{kneser},
and p. 13-26 in \cite{Pla}.
\index{cohomology}

\begin{dfn}
Let $\gal$ be a finite group,  $A$ a group provided with a $\gal$-action.
 $H^1(\gal,A)$ is defined as the quotient
\[
H^1(\gal,A)=\left\{\alpha:\gal\mapsto A\Big|\alpha(hg)=\alpha(h)
\alpha(g)^h\right\}/\equiv,
\]
where $\alpha\equiv\beta$ if and only if there exists
$a\in A$  such that $\alpha(g)=a^{-1}\beta(g)a^g$ for all $g\in\gal$.
If $\gal$ acts trivially on $A$, then $H^1(\gal,A)\cong\Hom(\gal,A)/A$,
where $A$ acts by conjugation.
In what follows we write $\alpha_g$ instead of $\alpha(g)$.
\end{dfn}

In case that $A\subseteq B$ is a subgroup, there is a long exact sequence
\[
 0\vaa A^\gal \vaa B^\gal \vaa (B/A)^\gal \vaa H^1(\gal,A) \vaa
 H^1(\gal,B),
\]
and furthermore, under the natural action
\footnote{ This result is not found in \cite{kneser}, but can be
found in \cite{Pla} p.22.}
of $B^\gal$ on $(B/A)^\gal$,
\begin{equation}\label{clasificador}
(B/A)^\gal/B^\gal \cong \ker\Big(H^1(\gal,A) \vaa H^1(\gal,B)\Big). \end{equation}

To simplify notations, in all that follows we assume that
whenever a sequence of pointed sets
\[
\dots\vaa U\vaa V\vaa W\vaa X\vaa Y \vaa Z
\]
is written, $X,Y,Z$ denote pointed sets,
$W,V,U,\dots$ denote groups,
and $W$ acts on $X$ with 
\[
X/W\cong\ker(Y\vaa Z).
\]
 
If $A$ is normal in $B$, we have in the sense just described
a long exact sequence
\begin{equation}\label{iso2}
\xymatrix{ 0\ar[r] & A^\gal\ar[r]  & B^\gal\ar[r]  & 
(B/A)^\gal \ar`r[d]`[lll]`[llld]`[llldr] [dll] \\
{}& H^1(\gal,A)\ar[r] & H^1(\gal,B)\ar[r] & H^1(\gal,B/A). }
\end{equation}
%

In case $A$ is central in $B$,
 the higher 
order cohomology groups for $A$ are also defined, and we have a long exact
 sequence
\[
\xymatrix{  0\ar[r]& A^\gal\ar[r]& B^\gal\ar[r]& 
(B/A)^\gal\ar`r[d]`[lll]`[llld]`[llldr] [dll]  \\
{} & H^1(\gal,A)\ar[r]& H^1(\gal,B)\ar[r]& 
H^1(\gal,B/A)\ar`r[d]`[lll]`[llld]`[llldr] [dll] \\
{}& H^2(\gal,A). &{}&{}}
\]

Finally, if $A$ and $B$ are both Abelian this sequence extends
 to cohomology 
of all orders (see \cite{Serre} or \cite{kenneth}).
All results of this section can be extended via direct limits
to profinite groups 
acting continuously on discrete groups
(see \cite{Serre}, p. 9 and p. 42).

We are specially interested in the case in which $\gal$ is 
the Galois group $\galkk$
of a  possibly infinite Galois extension $K/k$, where $k$ is a 
local or number field. In what follows, 
the subgroups $A,B,\dots$ 
are groups of algebraic or arithmetical nature.

The following  is a known fact, (see \cite{kneser},
ex. 1, p. 16).

\begin{prop}
\label{general}
For any finite dimensional algebra $A$, defined over $k$,
 and any algebraic extension
$K/k$, it holds that $H^1(\galkk,A^*_K)=\uno$.
\end{prop}

\begin{ex}[Hilbert's theorem 90]\label{Hilbert90}
 $\lineK \cong (\End_K(V))^*$.
Therefore,
 $H^1\Big(\gal,\lineK\Big)=\{1\}$.
\end{ex}

\section{Tensors and $K/k$-forms.}\label{kkforms}\index{tensors}

By a tensor of type $(l,m)$ on $V$,
we mean an $\Omega$-linear map
$\tau:V^{\otimes l}\vaa V^{\otimes m}$,
 where
\[
V^{\otimes r}=\bigotimes_{i=1}^rV
\textnormal{ for } r\geq1,\ V^{\otimes 0}=\Omega.
\]
$\tau$ is said to be defined over $k$, if $\tau(V_k^{\otimes l})\subseteq
V_k^{\otimes m}$. \emph{All tensors mentioned in this work are assumed to
be defined over $k$}.
 $\lineomega$ acts on the set of tensors of type
$(l,m)$ by $g(\tau)=g^{\otimes m}\circ\tau\circ(g^{\otimes l})^{-1}$.
It makes sense, therefore, to speak about the stabilizer of a tensor.

Let $I$ be any set. By an $I$-family of tensors, we mean a map that
associates, to each element $i\in I$, a tensor $t_i$ of type $(n_i,m_i)$.
$\lineomega$ acts on the set of all $I$-families by acting in each
coordinate.
In all that follows, we say a family instead of an $I$-family unless
the set of indices needs to be made explicit.
Let $\tn$ be a family of tensors and $H=Stab_{\lineomega}(\tn)$.
Then, $H$ is a linear algebraic group.

If $K/k$ is a Galois extension with Galois group $\gal$,
 we get an exact
 sequence 
\[
\uno\vaa H_K\vaa\lineK\vaa X_K\vaa\uno,
\]
where $X_K$ is the $\lineK$-orbit of $\tn$. 
It follows from (\ref{clasificador}),
and example \ref{Hilbert90}, that
$X_K^\gal/\linek\cong H^1(\gal,H_K)$.
The elements of  $X_K^\gal/\linek$ can be thought of as
isomorphism classes of pairs $(V'_k,\tn')$ that become isomorphic to
$(V_k,\tn)$ when extended to $K$. These classes are usually called
$K/k$-forms of $(V,\tn)$, or just $k$-forms if $K=\bark$.
\index{kkforms@$K/k$-forms}
Observe that two vector spaces of the same dimension are isomorphic,
so we can always assume that the vector space $V$, in which all 
tensors are defined, is fixed.

\begin{dfn}
We call a pair $(V,\tn)$, where $\tn$ is a family of tensors on $V$,
a \emph{space with tensors}, or simply a \emph{space}.
By abuse of language, we identify $(V,\tn)$ and $(V,\tn')$ whenever
$\tn,\tn'$ are in the same $\linek$-orbit,
i.e., if they correspond to the same $K/k$-form.
We say that $(V,\tn')$ is a $K/k$-form of $(V,\tn)$,
if $\tn$ and $\tn'$ are in the same $\lineK$ orbit.
\end{dfn}

\begin{ex}
Let $Q$ be a non-singular quadratic form on the space $V$. Then,
$\orthomega=Stab_{\lineomega}(Q)$.
Equivalence classes of non-singular
 quadratic forms on $V_k$ are classified by
$H^1\Big(\gal,\orthkbar\Big)$. A space, in this case, is what is usually called a
quadratic space.
\end{ex}

\section{Lattices.}\label{chapter5}

Let $k$ be a local or number field, $K/k$ a Galois extension,
 $G\subseteq\lineomega$
 an algebraic group
defined over $k$, $\la_k$ a lattice on $V_k$, $L_K$  a $\gal$-invariant
lattice on $V_K$. Let $\gal=\galkk$.

\begin{prop}\label{result1}
If there is an element $\varphi\in G_K$ such that
$\varphi(L_K)=\Lambda_K$, then
$a_\sigma=\varphi^\sigma\varphi^{-1}$ is a cocycle, and its class in
$H^1(\gal,G_K^\la)$ is independent of the choice of  $\varphi$, depending 
only on the orbit of $L_K$ under $G_k$.
The correspondence assigning, to every such
$G_k$-orbit of $\oink_K$-lattices,
 an equivalence class of cocycles,  is an injection.
The image of this map is
\[
\ker\Big(H^1(\gal,G_K^{\la})\stackrel{i_*}{\longrightarrow}H^1(\gal,G_K)\Big),
\]
where $i$ is the inclusion.
\end{prop}

\paragraph{Proof.}
$G_K$ acts on the set of $\oink_K$-lattices in $V_K$.
Let $X$ be the orbit of $\Lambda_K$.
 We have an exact sequence
\[
\uno\vaa G_K^{\la}\vaa G_K \vaa X\vaa\uno.
\]
Hence, by 
(\ref{clasificador}),
we get
$X^\gal/G_k \cong \ \ker\Big(H^1(\gal,G_K^{\la}) \vaa
H^1(\gal,G_K)\Big)$.
$\Box$

\begin{ex} Using the fact that $H^1\Big(\gal,\lineK\Big)=\{1\}$,
we obtain that the set of $\linek$-orbits of $\gal$-invariant $\oink_K$-lattices
isomorphic to $\la_K$ is in correspondence with $H^1\Big(\gal,\lineKlam\Big)$.
\end{ex}

If $G$ is defined as the stabilizer of a family of
 tensors, e.g., the unitary group of a hermitian form or the
automorphism group of an algebra, we get
a more precise result.

Recall that in section \ref{kkforms} we identified 
$K/k$-forms of $(V,\tn)$ with the corresponding 
$\linek$-orbits of families of tensors.

\begin{dfn}
Let $(V,\tn)$ be a space.
 A lattice in $(V_K,\tn)$
is a pair $(\la_K,\tn)$, where $\la_K$ is a lattice in $V_K$.
$\lineK$ acts on the set of pairs $(\la_K,\tn')$,
for all families of tensors $\tn'$, by acting on each component.
 Two lattices $(\la_K,\tn)$, $(L_K,\tn')$
are said to be \emph{in the same space} if $\tn,\tn'$
 are in the same $\linek$-orbit.
\end{dfn}

\begin{prop}\label{section3.4}
 Assume that $G$ is the
stabilizer of a family of tensors $\tn$ on $V$.
The set $H^1(\gal,G_K^{\la})$ is
in one-to-one correspondence
with the set of $G_k$-orbits of $\gal$-invariant
$\oink_K$-lattices in the same $G_K$-orbit,
in all spaces that are $K/k$-forms of $(V,\tn)$.
The kernel of the map
\[
H^1(\gal,G_K^\la)\stackrel{i_*}{\longrightarrow}H^1(\gal,G_K),
\]
where $i$ is the inclusion, corresponds to 
the subset of orbits of lattices that 
are in the
same space as $\la_K$.
\end{prop}

\paragraph{Proof.}
We have an action of $\lineK$ on the set of all pairs $(L_K,\tn')$,
where $L_K$ is a lattice and $\tn'$ an $I$-family of tensors with a
fixed index set $I$. If $T$ is the orbit of
$(\la_K,\tn)$, we have a sequence 
\[
\uno\vaa G_K^\la\vaa\lineK\vaa T\vaa\uno,
\]
and the same argument as before applies. Last statement
 follows from the fact that spaces $(V_K,\tn')$ are classified by
$H^1(\galkk,G_K)$, (see section \ref{kkforms} or \cite{kneser}, p. 15).
$\Box$

\begin{rmk}
Recall that $\la_K=\la_k\otimes_{\oink_k}\oink_K$.
If $L_K$ is in the same $G_k$-orbit as $\la_K$, $L_K=L_k
\otimes_{\oink_k}\oink_K$, since $G_k$ also acts on $V_k$.
Recall that we defined the cocycle corresponding to $L$
by the formula $a_\sigma=\varphi^\sigma\varphi^{-1}$ (see prop. 
\ref{result1}). This definition does not depend on $G$,
as long as $\varphi\in G$. 
It follows that the set of $G_k$-orbits of
lattices in $V_k$ that are isomorphic as $\oink_k$-modules,
and whose extensions to $K$ are in the same $G_K$ orbit,
corresponds to
\begin{equation}\label{kernelintersection}
\ker\Bigg(H^1(\gal,G_K^\la)\vaa H^1(\gal,G_K)\times
H^1\Big(\gal,\lineKlam\Big)\Bigg).
\end{equation}
In the case that $G$ is the stabilizer of a family of tensors,
\[
\ker\Bigg(H^1(\gal,G_K^\la)\vaa H^1\Big(\gal,\lineKlam\Big)\Bigg)
\]
corresponds to the set of $G_k$-orbits of
such lattices in all spaces that are
$K/k$-forms of $(V,\tn)$.
\end{rmk}

\begin{ex}
If $\la_k$ is free, (\ref{kernelintersection})
 corresponds to the set of $G_k$-orbits of free lattices
on $V_k$, whose extensions to $K$
are in the same $G_K$-orbit.
\end{ex}

\begin{dfn}
We say that an $\oink_K$-lattice $\la_K$ is 
\emph{defined over $k$}, if
$\la_K\cong\oink_K\otimes_{\oink_k}\la_k$ for some $\la_k$. 
We say that $\la_K$ is a \emph{$k$-free lattice}, if $\la_k$ is free.
\end{dfn}

Assume first that $G$ is the stabilizer of a family of tensors.

\begin{dfn}
Let $a\in H^1(\gal,\gKla)$. We say that $a$ is defined over 
$k$, $k$-free or
in $(V,\tn)$ if some (hence any),
lattice in the class corresponding to $a$ has this property. Define
\index{lk@$\lk,\lfree$, etc}
\[
\begin{array}{ccl}
\Lk G{K/k}{\la}&=&\{ a\in H^1(\gal,\gKla) | a \textnormal{
is defined over } k\},\\
\Lfree G{K/k}{\la}&=&\{a\in \Lk G{K/k}{\la}| a
\textnormal{  is $k$-free}\},\\
\Lvk G{K/k}{\la}&=&\{a \in  H^1(\gal,\gKla)| a
\textnormal{  is in } (V_K,\tn)\},\\
\Lkvk G{K/k}{\la}&=&\Lvk G{K/k}{\la}\cap \Lk G{K/k}{\la},\\
\Lfreevk G{K/k}{\la}&=&\Lvk G{K/k}{\la}\cap\Lfree G{K/k}{\la}.
\end{array}
\]
\end{dfn}

Let
\begin{equation}
\label{eq1}
F_1:  H^1(\gal,\gKla)\longrightarrow H^1(\gal,\gK),
\end{equation}
\begin{equation}
\label{eq2}
F_2:  H^1(\gal,\gKla)\longrightarrow H^1(\gal,\lineKlam),
\end{equation}
be the maps defined by the inclusions. Then, we have the following 
proposition:

\begin{prop}\label{proposition7}
Assume that $\la_k$ is free.
The following identities hold:
\[
\begin{array}{ccl}
\Lvk G{K/k}{\la} &=& \ker F_1,\\
\Lfree G{K/k}{\la} &=& \ker F_2,\\
\Lfreevk G{K/k}{\la} &=& \ker F_1\cap\ker F_2.\Box
\end{array}
\]
\end{prop}

Later, we give a similar interpretation to $\lk$.

\begin{ex}
\[
\Lfree {\orthomega}{\bark/k}{\la}=
\ker\Bigg(H^1(\gal,\orthi n{\bark}Q\la)\vaa H^1\Big(\gal,\linekbarlam\Big)\Bigg)
\]
is in correspondence with the set of isometry classes of
free quadratic lattices that become isometric to $\la_k$ over 
some extension.
\end{ex}

\begin{rmk}
Notice that $\lvk,\lkvk,\lfreevk$ can be defined, even if $G$
is not the stabilizer of a family of tensors, as follows:
\[
\begin{array}{rcl}
\Lvk G{K/k}{\la}&=&\ker\Big(H^1(\gal,\gKla)\vaa H^1(\gal,\gK)\Big),\\
\Lkvk G{K/k}{\la}&=&\Big\{a\in\Lvk G{K/k}{\la}\Big| a \textnormal{ is
 defined over }k\Big\} ,\\
\Lfreevk G{K/k}{\la}&=&\Big\{a\in\Lvk G{K/k}{\la}\Big| a \textnormal{ is
 free}\Big\}.
\end{array}
\]
In this case, the first and last identities of proposition
\ref{proposition7} still hold. Notice that we can still interpret
 $\lvk$ as a set of equivalence classes of lattices,
because of proposition \ref{result1}.
\end{rmk}

\paragraph{The set $H^1(\gal,U_K)$ and the ideal group.}\label{units}

Let $k$ be a local or number field,
$K/k$ a finite Galois extension. 
Let  $\gal=\galkk$, and let $U_K=\oink_K^*$ denote the group of units of $\oink_K$.
\index{uk@$U_k$}

For any local or number field $E$, let $I_E$ be its group of fractional
ideals, $P_E$ the subgroup of principal fractional ideals.
There is a natural map
$\alpha:I_k\rightarrow I_K$ defined by
$\alpha(\ideala)=\ideala\tensor\oink_K$.  Clearly
$\alpha(P_k)\subseteq P_K$,
so we get a map $\alpha': I_k/P_k\rightarrow I_K/P_K$.

We apply the general theory to $\la_k=\oink_k$, $G=\grupi_m$.  
Any $\lambda\in K^*$ acts by $\ideala\mapsto\lambda\ideala$, for $\ideala\in I_K$,
whence $G_K^{\la}=U_K$. It follows that,
\[
H^1(\gal,U_K)\cong(P_K)^\gal/\alpha(P_k).
\]

Non-zero prime ideals of $\oink_K$ form a set of free generators for $I_K$
(see \cite{Lang}, p. 18).
Let $\ideala\in I_K$. We can write 
\[
\ideala
=\prod_{\idealp\in \Pi(k)}(\prod_{\idealP | \idealp}
\idealP^{\beta(\idealP)}).
\]
 If $\ideala$ is $\gal$-invariant, all
the powers $\beta(\idealP)$ corresponding to prime
 divisors of the same prime of $k$ must be equal.  In other words:
\begin{equation}
\label{productus}
\ideala=\prod_{\idealp\in\Pi(k)}(\prod_{\idealP |
\idealp}\idealP)^{\beta(\idealp)},
\end{equation}
where $\beta(\idealp)$ is the common value of $\beta(\idealP)$ for all
$\idealP$ dividing $\idealp$.
This ideal is in $\alpha(I_k)$ if and only if the ramification degree
$e_\idealp$ divides
$\beta(\idealp)$ for all $\idealp$. Hence, we have an exact sequence
\[
0\longrightarrow \ker\alpha'\longrightarrow(P_K)^\gal/\alpha(P_k)
\longrightarrow\prod_{\idealp\in \Pi(k)}(\enteri/e_\idealp) ,
\]
where the image of the last map corresponds to those ideals of the form
 (\ref{productus}) that are principal in $K$.
The image of $\ker\alpha'$ is what we call
$\Lk G{K/k}{\la}$.
  In particular, since all ideals become principal in some
 extension, we can take a direct limit, to
obtain the long exact sequence:
\[
0\vaa I_k/P_k\vaa H^1(\galbarkk,U_{\bar{k}})
\vaa\prod_{\idealp\in \Pi(k)}(\Q/\enteri)\vaa 0.
\]
A refinement of this argument gives
\[
H^1(\galbarkk,U_{\bar{k}})\cong\left(I_k\otimes_\enteri\Q\right)/
\left(P_k\otimes_\enteri\enteri\right),\ \ \Lk G{K/k}{\la}=I_k/P_k.
\]

Recall remarks \ref{fixw} and \ref{galoisgroups}.
Assume $k$ is a number field.
There exist natural localization maps
\[
F_v:H^1(\gal,\gKla)\rightarrow H^1(\galv,\gKlav),
\]
defined by inclusion and restriction. 
We define $\gKlav=G_{K_w}$ if $w$ is Archimedean.

\begin{lem}
\label{localkernelgglam}
Let $F_1:H^1(\gal,\gKla)\rightarrow H^1(\gal,G_K)$ be the map induced by
the inclusion. If the natural map
\[
\tau:H^1(\gal,G_K)\rightarrow \prod_{v\in\Pi(k)} H^1(\galv,G_{K_w})
\]
is injective, then $\ker F_1\supseteq \bigcap_v \ker F_v$.
\end{lem}

\subparagraph{Proof of lemma.}
Immediate from the following
commutative diagram:
\[
\xymatrix{
H^1(\gal,\gKla) \ar[r]^{F_1}\ar[d]^{\prod_vF_v}&
H^1(\gal,G_K)\ar[d]^{\tau} \\
\prod_v H^1(\galv,\gKlav)\ar[r] &
\prod_v H^1(\galv,G_{K_w}).\Box}
\]

\begin{rmk}
If the hypothesis of this lemma is satisfied, one says that $G$
 satisfies the Hasse principle over $k$.
\end{rmk}

\subparagraph{Characterisation of $\lk$.}

$\Lk G{K/k}{\la}$ is the set of equivalence classes
 of lattices defined over $k$
that become isomorphic over $K$.
A lattice $L_K$ is defined over $k$ if and only if it is generated by its
$k$-points, i.e., 
\[
L_K=\oink_K(L_K\cap V_k).
\]
This is a local property, being an equality of lattices. On the other hand, 
for any local place $v$, all lattices defined over $k_v$ are $k_v$-free,
i.e., 
$$\Lk {\lineomega}{K_w/k_v}{\la}=\Lfree {\lineomega}{K_w/k_v}{\la}.$$
The following result is immediate from this observation.
\begin{prop}\label{lkernel}
\[
\Lk G{K/k}{\la}=\ker\left(H^1(\gal,\gKla)\vaa\prod_v
H^1\Big(\galv,\lineKvlam\Big)\right).\Box
\]
\end{prop}

\subsection{Genus and cohomology.}\label{chapter6}

Assume that
In all of section \ref{chapter6}, $k$ is a number field.

\begin{dfn}
Let $F_v$ be the localization map. Define\index{cgenus@$C$-genus}
\[
\ckgen G{K/k}{\la}=\ker\left(\prod_vF_v\right).
\]
We call this set the \emph{cohomological genus} of $\la$ with respect to
$G$. \index{cohomological genus}
\end{dfn}

\begin{prop}\label{wcisc}
For any linear algebraic group $G$, it holds that
\[
\ckgen G{K/k}{\la}\subseteq
\Lk G{K/k}{\la}.
\]
\end{prop}

\paragraph{Proof.}
This follows from proposition \ref{lkernel} and the 
commutative diagram
\[
\xymatrix{H^1(\gal,\gKla)\ar[d]\ar[rd]&\\
\prod_{v\in\Pi(k)}H^1(\galv,\gKlav)\ar[r]&
\prod_{v\in\Pi(k)}H^1\Big(\galv,\lineKvlam\Big).\Box}
\]

\begin{rmk}
Assume $G$ is the stabilizer of a family of tensors.
This result tells us that the cohomological genus corresponds to a 
set of equivalence classes ol lattices defined over $k$.
In fact, $a\in\ckgen G{K/k}{\la}$ if and only if $a$ corresponds to a lattice,
in some $K/k$-form of $(V,\tn)$, that is in the same $G_{k_v}$-orbit,
at every place $v$.
\end{rmk}

\begin{dfn}\index{vcgenus@$VC$-genus}
We define the $VC$-genus of $\la_k$
 by the formula
\[
V\ckgen G{K/k}{\la}=
\ckgen G{K/k}{\la}\cap \Lvk G{K/k}{\la}.
\]
In other words, it is the kernel of the map
\begin{equation}\label{cohomologymap}
H^1(\gal,\gKla)\vaa H^1(\gal.\gK)\times\prod_{v\in\Pi(k)}
H^1(\galv,\gKlav).
\end{equation}
\end{dfn}

Let $G$ be an arbitrary linear algebraic group.
The $VC$-genus corresponds to a set of $G_k$-orbits
of lattices in $V_k$. 
In fact, it corresponds to a subset of the set of double cosets
$\gk\backslash \gak/\gakla$, i.e., the \emph{genus} of $G$
(see \cite{Pla}, p. 440).
In particular, the following proposition holds.

\begin{prop}\label{prevstrongaproxy}
If $G$ has class number $1$ with respect to a lattice $\la_k$, then
(\ref{cohomologymap}) has trivial kernel for every Galois extension $K/k$
(compare with corollary 4 on p. 491 of \cite{Pla}).$\Box$
\end{prop}

This, in particular, applies to a group having absolute strong
approximation (see \cite{Pla}).
 However, we have a stronger result.

\begin{prop}
\label{strongaproxy}
If $G$ has absolute strong approximation over $k$, the map
(\ref{cohomologymap}) is injective.
\end{prop}

\subparagraph{Proof.}
Recall remark \ref{fixw}.
Let $M_K$,$L_K$ be two $\gal$-invariant $\oink_K$-lattices in
$V_K$, that are locally in the same $\gkv$-orbit for all $v$.
Then, we can choose
elements $\sigma_v\in\gkv$, such that $\sigma_v M_{K_w}=L_{K_w}$
for every place $v$, and $\sigma_v=1$ at all but a finite number of
 places. Now, any global element
$\sigma$, close enough to $\sigma_v$ at all finite
places where $\sigma_v\neq1$,
and stabilizing $M_{K_w}=L_{K_w}$ at the remaining finite places,
satisfies $\sigma M_K=L_K$, as claimed.
$\Box$

The following result is just a restatement of lemma
\ref{localkernelgglam}.

\begin{prop}\label{vcisc}
If $G$ satisfies the Hasse principle over $k$, then
\[
V\ckgen G{K/k}{\la}=\ckgen G{K/k}{\la}.\Box
\]
\end{prop}

This result tells us that, in the presence of Hasse principle, 
the cohomological genus corresponds to a subset of the genus
(compare with \cite{Rohlfs}, thm 3.3, p. 198).

 Let $G\subseteq GL(V)$ be a semi-simple group with universal cover
$\gtilde$ and fundamental group $\mu_n$. Let $K=\bark$.
The short exact sequence
\[
\uno\vaa \mu_n \vaa\gtilde_K\vaa G_K \vaa\uno,
\]
defines a map
$\theta:\gk\vaa H^1(\gal,F)=\ksobrekn$.

Let $\la_k$ be any lattice in $V_k$.
 The following proposition holds.

\begin{prop}
With the above notations,
$V\ckgen G{K/k}{\la}$ is in one-to-one correspondence with 
the genus of $G$ 
\textnormal{(compare with theorem 8.13 in \cite{Pla}, p. 490)}.
\end{prop}

\subparagraph{Proof.}
It suffices to show that any two $G_k$-orbits in the same genus are 
identified over some extension.
Without loss of generality, we assume $k$ is non-real. It suffices 
to check that they are in the same spinor genus (see \cite{spinor}).
Spinor genera are classified by 
\[
J_k/J_k^nk^*\Theta_{\ad}
(G^{\la}_{\ad_k}),
\]
where
 $\Theta_{\ad}(G^{\la}_{\ad_k})$
is the image of the local spinor norm
(see\footnote{ The case of an orthogonal group is already considered in \cite{Hsia99}.}
\cite{mithesis} or \cite{spinor}).
This is a finite set, and the representing adeles can be chosen to have
trivial coordinates at almost all places.
Therefore, it suffices to take an extension that contains the 
$n$-roots of unity, and $n$ roots of a finite set of local elements.$\Box$

This result allows us to use cohomology to study the genus of any
Semisimple group.

\section{Determinant class of a lattice.}\index{determinant class}

Let $[\ideala]$ be the $k^*$-orbit of the $\oink_K$-ideal $\ideala$.
Assume that 
\[
\la_k=\underbrace{\oink_k\oplus\dots\oplus\oink_k}_{n\textnormal{
\scriptsize times}}.
\]

The map
$\dete_*:H^1\Big(\gal,\lineKlam\Big)\vaa H^1(\gal,U_K)$ is the map induced
 in cohomology by the determinant.
It is surjective, since $\dete$ 
has a right inverse.
However, in general it is not injective, as the example below shows.

\begin{dfn}
 Let $L_K$ be a $\gal$-invariant lattice in $V_K$, and let $a$ be the
cocycle class corresponding to the $\linek$-orbit of $L_K$.  We
define the determinant class of $L_K$, which we denote $\dete_*(L_K)$, by:
\[
\dete_*(L_K)=\dete_*(a)\in H^1(\gal,U_K)\cong P_K^\gal/\alpha(P_k),
\]
and we identify it with the corresponding ideal class.
\end{dfn}

\begin{ex}
Using the standard embedding $GL(V)\times GL(W)\vaa GL(V\oplus W)$, 
it is easy to prove that  $\dete_*(\la_K\oplus L_K)=
\dete_*(\la_K)\dete_*(L_K)$. In particular, we obtain
that $\dete_*(\ideala_1\oplus\dots\oplus\ideala_n)=
[\ideala_1\dots\ideala_n]$.

Assume $k\subseteq K$ are local fields with maximal ideals $\idealp$,
$\idealP$. Assume that $\idealp\oink_K=\idealP^e$. Then,
\[
\dete_*(\underbrace{\idealP\oplus\dots\oplus\idealP}_e)=
[\idealP^e]=1=\dete_*(\underbrace{\oink_K\oplus\dots\oplus\oink_K}_e),
\]
but the latter lattice is defined over $k$
 and the first one is not. 
\end{ex}

Let $\lk=\Lk {GL(V)}{K/k}{\la}$.
We have the following result:

\begin{lem}
\label{tiroliro}
$\lk\cap\ker(\dete_*)=\uno$.
\end{lem}

\paragraph{Proof of lemma.}

This follows from the fact 
that all $k$-defined lattices are of the form 
$\ideala_k\oplus\oink_k\oplus\dots\oplus\oink_k$
(see \cite{Om}, (81:5)). 
It can also be proved by a diagram chasing argument.
$\Box$

Now observe that, for any algebraic group $G\subseteq GL(V)$,
we have $$\Lk G{K/k}{\la}=i^{-1}_*(\lk),$$
where $i_*$ is the cohomology map induced by the inclusion.
 
\begin{prop}
\label{lkernel2}
If $G\subseteq\specialomega$, then $i_*^{-1}(\lk)=\ker(i_*)$.
\end{prop}

\paragraph{Proof of proposition.}
It is immediate from the commutative diagram
$$\xymatrix{H^1(\gal,\gKla)\ar[rr]\ar[dr]_{i_*}&
{}\ar@{}[d]|{\circlearrowright}&H^1\Big(\gal,\specialKlam\Big)\ar[dl]
\\{}&H^1\Big(\gal,\lineKlam\Big)\ar[d]^{\dete_*}&{}\\&H^1(\gal,U_K)&{}}$$
 that $\im(i_*)\subseteq\ker(\dete_*)$.
Now recall  lemma \ref{tiroliro}.$\Box$

In particular, such a group cannot identify a free lattice to a
non-free $k$-defined lattice over any extension, although it can
identify a free lattice to a non-$k$-defined lattice.

In this case, a description of $\lfree$ is equivalent to a description
of $\lk$, hence $\Lk {G}{K/k}{\la}$ can be described without resorting to
localization.

\section{Example:Commutative algebras}

Let $A_k=k^n$ as a $k$-algebra. This is a space with a tensor of type $(2,1)$.
Then, the group of automorphisms of $A_K$ is the symmetric group $S_n$ for any
algebraic extension $K/k$. Asume henceforth that $K$ is an algebraic closure of
$k$.

The $k$-forms of $A$ are all semisimple commutative $k$-algebras.
This set is classified by $H^1(\gal,S_n)=Hom(\gal,S_n)/\equiv$, where $\equiv$
denotes the conjugation (as an equivalence relation).
Recall that any semisimple commutative algebra is a sum of fields.
If $\psi$ is the map that corresponds to an algebra $L_k=\bigoplus_iL_{i,k}$,
a simple computation shows that there is correspondence between
the fields $L_{i,k}$ and the orbits of $\im(\psi)$.

Let $R_k$ denote an $\oink_k$-subalgebra of $A_k$ (not necessarily with $1$)
 of maximal rank as a lattice.
Let $\Gamma=S_n^R$. Then, the image of the map $H^1(\gal,\Gamma)
\vaa H^1(\gal,S_n)$ corresponds to the set of isomorphisms classes of algebras
whose extensions to $K$ contain a $\gal$-invariant algebra isomorphic to $R_K$.
In particular, this includes all lattices defined over $k$, i.e.,
all $\oink_k$-algebras whose extensions to $K$ are isomorphic to $R_K$.

\begin{ex}
If $\Gamma$ is not transitive, then no field can contain an $\oink_k$-algebra
whose extension to $K$ is isomorphic to $R_K$. This is the case for example
of the algebra $R_k=\oink_k^{n-1}\oplus\ideala$ or
$R_k=\{a\in \oink_k^n|a_{n-1}-a_n\in\ideala\}$ if $\ideala$ is an ideal different
from $(1)$, and $n\geq3$.
\end{ex}

\begin{rmk}
All result in this paper apply also to lattices over rings of 
$S$-integers. Absolute strong approximation must be replaced by strong
 approximation with respect to $S$.
\end{rmk}

\printindex

\end{document}